\documentclass[12pt,draft]{article}

\begin{document}

\begin{center}
\LARGE\noindent\textbf{On Hamiltonian Bypasses in Digraphs satisfying Meyniel-like Condition}\\

\end{center} 

\begin{center}
\noindent\textbf{Samvel Kh. Darbinyan}\\

Institute for Informatics and Automation Problems of NAS RA
 
E-mail: samdarbin@iiap.sci.am \\

A translation from Russian of the paper by Darbinyan (Mathematical Problems of Computer Science,  vol. 20,  7-19, 1998) with some additional information.

\end{center}

\textbf{Abstract}

Let $G$ be a strongly connected directed graph of order $p\geq 3$. In this paper, we show that if $d(x)+d(y)\geq 2p-2$ 
(respectively, $d(x)+d(y)\geq 2p-1$) for every pair of non-adjacent vertices $x, y$, then $G$ contains a Hamiltonian path   (with only a few  exceptional cases that can be clearly   characterized) in which the initial vertex dominates the terminal vertex    (respectively, $G$ contains two distinct verteces $x$ and $y$ such that there are two internally disjoint $(x,y)$-paths of lengths $p-2$ and $2$).


\section {Introduction} 

In this paper we consider finite digraphs (directed graphs) without loops and multiple arcs. Every cycle and path is assumed simple and directed. We shall assume that the reader is familiar with the the standard terminology on digraphs  and refer to  \cite{[11]}  for  terminology and notation not described  in this paper.

A digraph $G$ of order $p$ is {\it Hamiltonian} (respectively, is {\it pancyclic})  if it contains a {\it Hamiltonian cycle},
i.e., a cycle that includes every vertex of $G$ (respectively, contains cycles of all lengths $m$, $3\leq m\leq p$).\\
 
\noindent\textbf{Definition 1.1.} { \it Let $G$ be a digraph of order $p$, and let $k$ be an integer. We will  say that a digraph $G$  satisfies condition $(M_k)$ if 
 $$d(x)+d(y)\geq 2p-2+k$$ 
 for every pair of non-adjacent vertices $x$, $y$ of $G$.}\\

Many researchers investigated hamiltonicity and pancyclcity of digraphs with condition $(M_k)$, $k\geq 0$ (see, e.g., \cite{[4]}, \cite{[6]}). 
  We now recall the following well-known degree conditions (Theorems 1.2 - 1.5) that guarantee that a digraph is Hamiltonian.\\ 

 \textbf{Theorem 1.2} (Nash-Williams \cite{[13]}). {\it Let $G$ be a digraph of order $p$ such that for every vertex $x\in V(G)$, $od(x)\geq p/2$ and $id(x)\geq p/2$. Then $G$ is Hamiltonian.}

 \textbf{Theorem 1.3} (Ghouila-Houri \cite{[9]}). {\it Let $G$ be a strong  digraph of order $p$ such that for every $x\in V(G)$, $d(x)\geq p$.  Then $G$ is Hamiltonian.}

Note that Theorem 1.2 is a consequence of Theorem 1.3.\\

 \textbf{Theorem 1.4} (Woodall \cite{[16]}). {\it Let $G$ be a digraph of order $p\geq 2$. If $od(x)+id(y)\geq p$ for all pairs of vertices $x$ and $y$ such that there is no arc from $x$ to $y$, then $G$ is Hamiltonian.}

\textbf{Theorem 1.5} (Meyniel \cite{[12]}). {\it Let $G$ be a strong digraph of order $p\geq 2$ satisfying condition $(M_1)$.  Then $G$ is Hamiltonian.}\\

Note that Meyniel's theorem is a generalization of Nash-Williams', Ghouila-Houri's and Woodall's theorems. For a short proof of Theorem 1.5, see \cite{[5]}. 
Nash-Williams \cite{[13]} raised the following problem.

\textbf{Problem 1.6} (Nash-Williams \cite{[13]}). {\it Describe all the extreme digraphs for the Ghouila-Houri theorem, i.e., describe all the strong non-Hamiltonian digraphs of order $p$ with minimum degree $p-1$}.\\

As a partial solution to Problem 1.6,  Thomassen proved a structural theorem on the extremal digraphs.

\textbf{Theorem 1.7} (Thomassen \cite{[14]}). {\it Let $G$ be a strong non-Hamiltonian digraph of order $p\geq 3$ with  minimum degree  $p-1$.  Let $C_m=x_1x_2\ldots x_mx_1$ be a longest cycle in $G$. 
Then 
any two distinct vertices of $V(G)\setminus V(C_m)$ are adjacent, every vertex of $V(G)\setminus V(C_m)$ has degree  $p-1$ in $G$, and 
every strong component of $G\langle V(G)\setminus V(C_m) \rangle $ is a complete digraph. Furthermore, if $G$ is 2-strong, then $C_m$ can be chosen such that $G\langle V(G)\setminus V(C_m) \rangle$ is a transitive tournament.}\\

It is natural to consider the analogous problem for the Meyniel theorem. In \cite{[7]}, we proved  Theorem 1.9. The following notation will be used in Theorem 1.9.

\textbf{Notation 1.8}. {\it For any $k\in [1, p-2]$ let $D_{p-k, k}$ denote a digraph of order $p\geq 3$, obtained from $K^*_{p-k}$ and $K^*_{k+1}$ by identifying a vertex of the first with a vertex of the second.}

\noindent\textbf{Theorem 1.9} (Darbinyan \cite{[7]}, for a detailed proof, see, arXiv:1911.05998v1). {\it Let $G$ be a strong non-Hamiltonian digraph of order $p\geq 3$ satisfying condition $(M_0)$. Let $C_m=x_1x_2\ldots x_mx_1$ be a longest cycle in $G$ and let $G_1$, $G_2, \ldots , G_h$ be the strong components of
$G\langle V(G)\setminus V(C_m)\rangle$ labelled in such a way that no vertex of $G_i$ dominates a vertex of $G_j$ whenever $i>j$.
Then the following statements hold:

I. Any two distinct vertices of $A:=V(G)\setminus V(C_m)$ are adjacent; every vertex of $A$ has degree at most $p-1$ in $G$; and 
every component $G_i$ $($$1\leq i\leq h$$)$  is a complete digraph.

II. If $G$ is not isomorphic to $D_{p-k, k}$, where $k\in [1, p-2]$, then for every $l\in [1,h]$ there are two distinct vertices $x_i,x_j$ on
 $C_m$  and  some vertices $u,v$ in $V(G_l)$ $($possibly, $u=v$$)$ such that $x_iu, vx_j\in E(G)$ (by $B_l$ is denoted the set $V(C_m[x_{i+1},  x_{j-1}])$)  and   
$$
E(B_l\rightarrow V(G_1)\cup V(G_2)\cup \cdots \cup V(G_l)) 
=E(V(G_l)\cup V(G_{l+1})\cup \cdots \cup V(G_h)\rightarrow B_l)=\emptyset,
$$
in particular, $E(V(G_l),B_l)=\emptyset$.  Moreover, $G\langle B_l\rangle$ also is a complete digraph;  
$$
V(G_1)\cup V(G_{2})\cup \cdots \cup V(G_{l-1})\rightarrow  B_l\cup V(G_l)\rightarrow
V(G_{l+1})\cup \cdots \cup V(G_h);
$$
and for all vertices $z\in V(G_l)$ and $y\in B_l$, $d(z,V(C_m))=m-|B_l|+1$ and  $d(y,V(C_m))=m+|B_l|-1$.
Any vertex of $B_l\cup V(G_l)$ cannot be inserted into $C_m[x_j,x_i]$ $($in particular, 
$x_i\rightarrow B_l\cup V(G_l)\rightarrow x_j$$)$.

III. If $D$ is 2-strong, then $G\langle A\rangle$ is a transitive tournament.

IV. For every  $r\in [2,m]$, $G$ contains cycle of  length $r$, unless when $p$ is odd  and $G$ is isomorphic to the complete bipartite digraph $ K^*_{\lfloor p/2\rfloor,\lfloor p/2\rfloor+1}$}. \\

Observe that, Theorem 1.7 follows directly from Theorems 1.9(I), 1.9(III).
Thomassen \cite{[15]} investigated the pancyclicity of digraphs with condition $(M_2)$, and  proved the following theorem.

\textbf{Theorem 1.10} (Thomassen \cite{[15]}). Let $G$ be a strong digraph of order $p\geq 3$ satisfying condition $(M_2)$. Then $G$ is pancyclic unless  $p$ is even and $G$ is isomorphic to $K^*_{p/2,p/2}$.\\

In \cite{[6]}, the author characterized those strong digraphs which satisfy Meyniel's condition (condition $(M_1)$), but are not pancyclic. Before stating the main result of \cite{[6]},  we need to define a family
 $\Phi^m_p$ of digraphs (see also in \cite{[15]}). 

\textbf{Notation 1.11}. {\it For any $p\geq 4$ and $m$, where $(p+1)/2<m\leq p-1$, by $\Phi^m_p$ we denote the set of digraphs $G$ with vertex set $\{x_1,x_2,\ldots , x_p\}$, which satisfy the following conditions:

a. $G$ satisfies condition $(M_1)$;

b. $x_1x_p\in E(G)$ and $x_{i+1}x_i\in E(G)$ for all $i\in [1,p-1]$ $($i.e., $x_px_{p-1}\ldots x_2x_1x_p$ is a Hamiltonian cycle in $G$$)$;

c. $E(x_i,x_{i+m-1})=\emptyset$ and $x_jx_i\notin E(G)$ whenever $2\leq  i+1<j\leq p$.}\\

Note that if $D\in \Phi^m_p$, then $D$ contains a cycle of length $k$, for every $k\in [2,p]\setminus \{m\}$, and it contain no cycly of length $m$.\\

\textbf{Theorem 1.12} (Darbinyan \cite{[6]}). Let $G$ be a strong digraph of order $p\geq 3$ satisfying condition ($M_1$). Then one of the following holds:

i. $G$ is pancyclic;

ii. $p$ is even and $G$ is isomorphic to $K^*_{p/2,p/2}$ or
$K^*_{p/2,p/2}\setminus \{u\}$, where  $u$ is an arbitrary arc of $K^*_{p/2,p/2}$; 

iii. $G\in \Phi^m_p$ for some $m$, $(p+1)/2<m\leq p-1$.\\

Later on, Theorem 1.12 also was proved by Benhocine \cite{[2]}. In \cite{[2]}, the author mentions that in \cite{[3]}, it was showed: If a digraph $D$ of order $p$ satisfies condition $(M_2)$, then it contains bypasses of every length $n$, $3\leq n\leq p$, with some exception (A bypass of length $n$ is a digraph obtained from a directed  cycle of length $n$ by reversing exactly one arc).

For any $n$ and $p$, where $2\leq n\leq p-2$, by $D(p,n)=[x_1x_2\ldots x_n; x_1y_1y_2\ldots$ $ y_{p-n}x_n]$ we denote a digraph of order $p$ with vertex set $\{x_1,x_2,\ldots , x_n,y_1,y_2,\ldots , y_{p-n}\}$ and arc set 
$$
\{x_1y_1, y_{p-n}x_n\}\cup \{x_i x_{i+1}\, |\, 1\leq i\leq n-1\} \cup \{y_iy_{i+1} \, |\, 1\leq i\leq p-n-1\}.
$$
In particular, $D(p,2)$ (respectively, $D(p,3)$) is a digraph obtained from a cycle $C$ of length $p\geq 3$ by reversing 
exactly one arc (respectively, exactly two consecutive arcs).
If a digraph $G$ of order $p$ contains a $D(p,2)$, then $D(p,2)$ is called a  Hamiltonian bypass in $G$. \\

It is natural to consider the following problem:

\textbf{Problem 1.13}. Whether a digraph  $D$ of order $p$  contains a Hamiltonian bypass (or $D$ contains a $D(p,n)$, where $3\leq n\leq p-2$)?\\

Benhocine \cite{[1]} proved that if a digraph $G$ satisfies the conditions of Theorem 1.2 or 1.3 or 1.4, 
 then $G$ contains a Hamiltonian bypass.\\ 

In view of the next theorems we need the following notations.

\textbf{Notation 1.14}. {\it Let $D_0$ denotes any digraph of order $p\geq 3$, $p$ is odd, such that $V(D_0)=A\cup B$, where $A\cap B=\emptyset$, $A$ is an independent set with $(p+1)/2$ vertices, $B$ is a set of $(p-1)/2$ vertices inducing an arbitrary sudigraph, and $D_0$ contains all the possible arcs between $A$ and $B$.}

\textbf{Notation 1.15}. {\it By $T_5$ we denote a tournament of order 5 with vertex set $\{z_1, z_2,$ $ z_3, z_4, y\}$ and arc set $\{z_iz_{i+1} \, |\, i\in [1,3]\}\cup \{z_4z_1, z_1y, z_3y, yz_2, yz_4,z_1z_3, z_2z_4\}$.}\\

 The tournament $T_5$ satisfies condition $(M_0)$, but has no Hamiltonian bypass. Notice that $T_5$ also is not 2-strong, as $id(z_1)=od(z_4)=1$.\\

Benhocine \cite{[1]}  also proved the following theorem:

\textbf{Theorem 1.16} (Benhocine \cite{[1]}). {\it Let $G$ be a  2-strong digraph of order $p$  with minimum degree at least $p-1$. Then $G$ contains a Hamiltonian bypass, unless $G$ is isomorphic to a digraph of type $D_0$.}\\

An {\it oriented graph} is a digraph with no cycle of length two. In  \cite{[8]}, we  studied the problem of the existence of $D(p,3)$ in  oriented graphs with the large in-degrees and out-degrees. We proved the following theorem.

\textbf{Theorem 1.17} (Darbinyan \cite{[8]}). {\it Let $G$ be a oriented graph of order $p\geq 10$. If the minimum in-degree and out-degree of $G$ at least $(p-3)/2$, then $D$ contains a $D(p,3)$.}\\

In \cite{[1]}, Benhocine notece that to prove the existence of $D(p,2)$ in digraphs satisfying the degree conditions of the Meyniel theorem, seeme defficult to extend. 
In this paper we prove the following two theorems.

\textbf{Theorem 1.18.} {\it Let $G$ be a strong digraph of order $p\geq 3$ satisfying condition $(M_0)$.  Then $D$ contains a Hamiltonian bypass unless $G$ is isomorphic to a digraph of type  $D_0$ or $D_{p-k,k}$ or $G\in \{T_5, C_3\}$.}\\

  \textbf{Theorem 1.19.}  {\it Let $G$ be a strong  digraph of order $p\geq 4$ satisfying condition $(M_1)$. 
 Then $G$ contains a $D(p,3)$.}\\
 
Since none of the digraphs $T_5$, $C_3$ and $D_{p-k,k}$ are not 2-strongly connected,  Theorem 1.16  is an immediate consequence of  Theorem 1.18.
Using Theorem 1.18, it is not difficult to prove that if a digraph $G$ satisfies condition $(M_1)$, then it contains a Hamiltonian bypass, unless when $G\in \{C_3, T_5\}$ (see Corollary 3 in Section 4).
The  last two results (Theorems 1.18 and 1.19) also were presented at 5-th Scienece-Technical Conferense, Tsaghkadzor, Armenia, 1986.

\section {Further terminology and notations}

  In this paper we consider finite digraphs without loops and multiple arcs. For a digraph $G$, we denote
  by $V(G)$ the vertex set of $G$ and by  $E(G)$ the set of arcs in $G$. The {\it order} of $G$ is the number
  of its vertices. Let $x$, $y$ be distinct vertices in $G$.
 The arc of a digraph $G$ directed from
   $x$ to $y$ is denoted by $xy$. For disjoint subsets $A$ and  $B$ of $V(G)$  we define $E(A\rightarrow B)$ 
   as the set $\{xy\in E(G)\, |\, x\in A, y\in B\}$.   
If $x\in V(G)$
   and $A=\{x\}$ we write $x$ instead of $\{x\}$. 
 The {\it out-neighborhood} of a vertex $x$ is the set $O(x)=\{y\in V(G)\, |\, xy\in E(G)\}$ and $I(x)=\{y\in V(G)\, |\, yx\in E(G)\}$ is the {\it in-neighborhood} of $x$. Similarly, if $A\subseteq V(G)$, then $O(x,A)=\{y\in A \,| \, xy\in E(G)\}$ and $I(x,A)=\{y\in A \, |\, yx\in E(G)\}$. 
The {\it out-degree} of $x$ is $od(x)=|O(x)|$ and $id(x)=|I(x)|$ is the {\it in-degree} of $x$. Similarly, $od(x,A)=|O(x,A)|$ and $id(x,A)=|I(x,A)|$. The {\it degree} of the vertex $x$ in $G$ defined as $d(x)=od(x)+id(x)$ (similarly, $d(x,A)=od(x,A)+id(x,A)$). 

The subdigraph of $G$ induced by a subset $A$ of $V(G)$ is denoted by $G\langle A\rangle$. For integers $a$ and $b$, $a\leq b$, by  $[a,b]$  we denote  the set $\{a,a+1,\ldots , b\}$. 
The  path (respectively, the  cycle) consisting of the distinct vertices $x_1,x_2,\ldots ,x_m$ ($m\geq 2 $) and the arcs $x_ix_{i+1}$, $i\in [1,m-1]$  (respectively, $x_ix_{i+1}$, $i\in [1,m-1]$, and $x_mx_1$), is denoted  $x_1x_2\cdots x_m$ (respectively, $x_1x_2\cdots x_mx_1$). We say that $x_1x_2\cdots x_m$ is a {\it path from $x_1$ to $x_m$} or is an  $(x_1,x_m)$-{\it path}. The {\it length} of a cycle or a path is the number of its arcs. A cycle of length $k$, $k\geq 2$, is denoted by $C_k$. For a cycle  $C_k:=x_1x_2\cdots x_kx_1$, the subscripts considered modulo $k$, i.e. $x_i=x_s$ for every $s$ and $i$ such that  $i\equiv s\, (\hbox {mod} \,k)$. 

A cycle (path) that contains  all the vertices of a digraph $G$  is a {\it Hamiltonian cycle} (is a {\it Hamiltonian path}).   
 A digraph $G$ is {\it strongly connected} (or, just, {\it strong}) if there exists a path from $x$ to $y$ and a path from $y$ to $x$ for every pair of distinct vertices $x,y$. A digraph $G$ is $k$-{\it strongly connected} (or, $k$-{\it strong}), if $|V(G)|\geq k+1$ and  $G\langle V(G)\setminus A\rangle$ is strong for any set $A$ of at most $k-1$ vertices. 

  For an undirected graph $G$, we denote by $G^*$ symmetric digraph obtained from $G$ by replacing every edge $xy$ with the pair $xy$, $yx$ of arcs.  $K_{p,q}$ denotes the complete undirected bipartite graph  with partite sets of cardinalities $p$ and $q$.  Two distinct vertices $x$ and $y$ in a digraph $G$ are {\it adjacent} if $xy\in E(G)$ or $yx\in E(G) $ (or both).

\section { Preliminaries }

The following well-known simple Lemmas 3.1-3.4 are the basis of our results
and other theorems on directed cycles and paths in digraphs. They
will be used extensively in the proof of our result. 

\textbf{Lemma 3.1} (H\"{a}ggkvist and Thomassen \cite{[10]}). {\it Let $G$ be a digraph of order  $p\geq 3$
 containing a
 cycle $C_m$, $m\in [2,p-1] $. Let $x$ be a vertex not contained in this cycle. If  $d(x,V(C_m))\geq m+1$,
 then  for every  $k\in [2,m+1]$, $G$ contains a cycle of length $k$ including $x$.}\\

The following lemma is a  modification of a lemma by Bondy and Thomassen \cite{[5]}. 

\textbf{Lemma 3.2}. {\it Let $G$ be a digraph of order $p\geq 3$ containing a
 path $P:=x_1x_2\ldots x_m$, $m\in [2,p-1]$ and  $x$ be a vertex not contained in this path.
  If one of the following conditions holds:

 $(i)$ $d(x,V(P))\geq m+2$;

 $(ii)$ $d(x,V(P))\geq m+1$ and $xx_1\notin E(G)$ or $x_mx\notin E(G)$;

 $(iii)$ $d(x,V(P))\geq m$, $xx_1\notin E(G)$ and $x_mx\notin E(G)$};

\noindent\textbf{} {\it then there is an  $i\in [1,m-1]$ such that
$x_ix,xx_{i+1}\in E(D)$  i.e.,  $x_1x_2\ldots
x_ixx_{i+1}\ldots x_m$ is a path of length $m$ in $G$ $($we say that  $x$ can be
inserted into $P$ or the path $x_1x_2\ldots x_ixx_{i+1}\ldots$ $ x_m$
is an  extended path obtained from $P$ with $x$$)$.}\\

 Using Lemma 3.1 (respectively, Lemma 3.2) one can prove  Lemma 3.3 (respectively,  Lemma 3.4).

\textbf{Lemma 3.3}. {\it Let $G$ be a digraph of order $p\geq 3$ and   $C_m$ be a cycle in $G$, where $2\leq m\leq p-1$. Suppose that for every vertex $y\in V(G)\setminus V(C_m)$, $d(y)\geq 2p-m-1$.  Then for any subset $A\subseteq V(G)\setminus V(C_m)$, $G$ contains a cycle with vertex set $A\cup V(C_m)$.}

\textbf{Lemma 3.4}. {\it  Let $G$ be a digraph of order $p\geq 3$. Suppose that $P:=x_1x_2\ldots x_m$ is a path in $G$, where $2\leq m\leq p-1$, and for every vertex $y\in V(G)\setminus V(P)$, $d(y)\geq 2p-m$.  Then for any subset $A\subseteq V(G)\setminus V(P)$, $G$ contains an  ($x_1, x_m$)-path with vertex set $A\cup V(P)$.}\\

In the proof of  Theorem 1.18 we also need the following
 lemma. 

\textbf{Lemma 3.5}. {\it Let $G$ be a digraph of order $p\geq 3$ and let $C:=x_1x_2\ldots x_{p-1}x_1$ be a cycle of length $p-1$ in $G$. Suppose that $y\notin V(C)$  and $G$ contains no Hamiltonian bypass. Then the following statements hold:

$(i)$ $od(y,\{x_i, x_{i+1}\})\leq 1$ and   $id(y,\{x_i, x_{i+1}\})\leq 1$ for all $i\in [1,p-1]$;

$(ii)$ $od(y)\leq (p-1)/2$, $id(y)\leq (p-1)/2$ and $d(y)\leq p-1$;

$(iii)$ if $k\in [1,p-1]$ and  $x_ky, yx_{k+1}\in E(G)$, then $x_{i+1}x_i\notin E(G)$ for all $i\in [1,p-1]\setminus \{k\}$.}\\

\section {Proofs of the main results}

The following definition will be used in our proofs.

\textbf{Definition 4.1}. {\it Let $P_0:=x_1x_2\ldots x_m$, $m\geq 2$, be an $(x_1,x_m)$-path in a digraph $G$. Assume that  the vertices
$y_1,y_2,\ldots , y_k$ are in $V(G)\setminus  V(P_0)$, $k\geq 1$. For $i\in [1,k]$, by $P_i$ we denote  an $(x_1,x_m)$-path in $G$ with vertex set $V(P_{i-1})\cup \{y_j\}$ $($if it exists$)$, i.e., $P_i$ is an extended path obtained from $P_{i-1}$ with some vertex  $y_j\notin V(P_{i-1})$. If $e+1$ is the maximum possible number of these paths $P_0, P_1,\ldots , P_e$, $e\in [0,k]$, then we say that $P_e$ is an extended path obtained from $P_0$ with vertices $y_1,y_2,\ldots , y_k$ is as much as possible. Notice that for all $i\in [0,e]$}, $P_i$ is an $(x_1,x_m)$-path of length $m+i-1$.\\

\textbf{Proof of  Theorem 1.18}. 

Let $G$ be a digraph of  order $p\geq 3$. It is clear that one of the following holds:

i. $G$ contains a cycle of length $p-1$;

ii. $G$ is Hamiltonian and contains no cycle of length $p-1$;

iii. The longest cycles in $G$ has length at most $p-2$.\\

Because of these, to prove Theorem 1.18 it suffices to prove the following Lemmas 4.2-4.4 below.

\textbf{Lemma 4.2}. {\it Let $G$ be a digraph  of order $p\geq 3$ satisfying
 condition $(M_0)$. If $G$  contains a cycle of length $p-1$, then either $G$ contains a Hamiltonian bypass or $G\in D_0\cup \{D_{p-1,1},T_5\}$}.

\textbf{Lemma 4.3}. {\it Let $G$ be a Hamiltonian digraph of order  $p\geq 3$ satisfying condition $(M_0)$. If $G$  contains no cycle of length $p-1$, then either $G$ contains a Hamiltonian bypass or it is isomorphic to the directed cycle of length three}.

\textbf{Lemma 4.4}. {\it Let $G$ be a strong non-Hamiltonian digraph of order $p\geq 3$ satisfying condition $(M_0)$. 
If $G$ contains no cycle of length $p-1$, then  either $G$  contains  a Hamiltonian bypass or it is isomorphic to the digraph 
$D_{p-k,k}$, where $1\leq k\leq p-2$}.\\ 

\textbf{Proof of Lemma 4.2}. 

 Let $G$ be a digraph of order $p\geq 3$ satisfying the conditions of Lemma 4.2. Suppose  that  $G$ contains no Hamiltonian bypass and  $G\notin D_0\cup \{D_{p-1,1},T_5\}$. 
Let $C:=x_1x_2\ldots x_{p-1}x_1$ be an arbitrary cycle of length $p-1$ in $G$ and  $y$ be the vertex  not in $C$. For the cycle $C$ and the vertex $y$, let us prove Claims 1 and 2.

\textbf{Claim 1.} There are no pair of integers $s\in [1,p-1]$ and $k\in [2,p-2]$ such that 
$$
x_sy,yx_{s+k}\in E(G) \quad \hbox{and} \quad 
E(y,\{x_{s+1},x_{s+2}, \ldots ,  x_{s+k-1}\})=\emptyset.  
$$
{\it Proof}. By contradiction, suppose  that there exist  some $s\in [1,p-1]$ and $k\in [2,p-2]$ such that the following holds:
 $$
x_sy,yx_{s+k}\in E(G) \quad \hbox{and} \quad 
E(y,\{x_{s+1},x_{s+2}, \ldots ,  x_{s+k-1}\})=\emptyset.   \eqno (1)
$$
Choose $s$ and $k$ such that $|\{x_{s+1},x_{s+2}, \ldots ,  x_{s+k-1}\}|$ be the smallest possible. We will consider the cases $k=2$ and $k\geq 3$ separately.

\textbf{Case 1.} $k=2$.

From $E(y,x_{s+1})=\emptyset$, Lemma 3.5(ii) and condition $(M_0)$ it follows that $d(y)=d(x_{s+1})=p-1$. Now using Lemma 3.5(i),
we obtain $od(y)= id(y)= (p-1)/2$. Therefore,  $p$ is odd and
$$
O(y)=I(y)=\{x_{s},x_{s+2}, \ldots ,  x_{s-2}\}. \eqno (2)
$$
Let $x_j$ be an arbitrary vertex in $\{x_{s+1},x_{s+3}, \ldots ,  x_{s-1}\}$. From condition ($M_0$), (2) and Lemma 3.5(ii) it follows that $d(x_j)=p-1$ and $x_{j-1}yx_{j+1}x_{j+2}\ldots x_{j-2}x_{j-1}$ is a cycle of length $p-1$. Similarly to (2), we can show 
that 
$
O(x_j)=I(x_j)=\{x_{j+1},x_{j+3}, \ldots ,  x_{j-1}\} 
$.
Therefore, $\{x_{s+1},x_{s+3}, \ldots ,  x_{s-1},y\}$ is an independent set, i.e., $G$ is isomorphic to a digraph of type  $D_0$.

\textbf{Case2.} $k\geq 3$.

We first show that 
$$
d(y)\leq p-k. \eqno (3)
$$
Assume that (3) is false, i.e., $d(y)\geq p-k+1$. Using (1) and Lemma 3.5(i) it is not difficult to show that $p-k$ is odd, 
 $od(y)=id(y)=(p-k+1)/2$ and
$$
O(y)=I(y)=\{x_{s+k},x_{s+k+2}, \ldots ,  x_{s-2},x_s\}. 
$$
Therefore, $x_{s+k}y, yx_{s+k+2}\in E(G)$ and $E(y,x_{s+k+1})=\emptyset$, which contradict that $k$ is minimal. Therefore,
$d(y)\leq p-k$.

Now from (1), (3) and condition ($M_0$) it follows that $d(x_{s+l})\geq p+k-2$ for all $l\in [1,k-1]$. Observe that $Q:=x_syx_{s+k}\ldots  x_{s-1}x_s$ is a cycle of length $p-k+1$. Since 
$$
 2p-|V(Q)|-1=p+k-2\leq d(x_{s+l}),
$$ 
we can apply Lemma 3.3 to the cycle $x_syx_{s+k}x_{s+k+1}\ldots x_{s-1}x_{s}$ and the set $\{x_{s+1},x_{s+2},$  $\ldots , x_{s+k-1}\}$.
In a result, we obtain a cycle of length $p-1$, which does not contain the vertex $x_{s+k-1}$. This contradicts Lemma 3.5(ii), since
$D$ contains no Hamiltonian bypass and $d(x_{s+k-1})\geq p+k-2\geq p$. Claim 1 is proved. \fbox \\

\textbf{Claim 2}. $d(y)=p-1$.

{\it Proof}. By contradiction, suppose   that  $d(y)\not=p-1$. This together with Lemma 3.5(ii) implies that $d(y)\leq p-2$. Therefore, there are integers $s\in [1,p-1]$ and $k\geq 2$ such that 
$$
E(y,\{x_{s+1},x_{s+2},\ldots , x_{s+k-1}\})=\emptyset; \eqno (4)
$$
$$
E(y,x_{s})\not=\emptyset \quad \hbox{and}\quad E(y,x_{s+k})\not=\emptyset,  \eqno (5)
$$
where $k\geq 2$. Since $G\not\cong D_{p-1,1}$, we have $x_s\not= x_{s+k}$ and $p\geq 4$. By Claim 1,
$$
|E(x_s\rightarrow y)|+|E(y\rightarrow x_{s+k})|\leq 1. \eqno (6)
$$
It is not difficult to show that
$$ d(y)\leq p-k.   \eqno (7)   $$
Indeed, if $p-k$ is even, then (7) immediately follows from Lemma 3.5(i). We may therefore assume that $p-k$ is odd. Then using (5) and (6), we obtain
either $x_sy\in E(G)$ and $yx_{s+k}\notin E(G)$, or $x_sy\notin E(G)$. Thus, we have that  $yx_{s+k}\notin E(G)$ or $x_sy\notin E(G)$. Now, again using Lemma 3.5(i), it is easy to see that: 

If $x_sy\notin E(G)$, then $id(y)\leq (p-k-1)/2$ and $od(y)\leq (p-k+1)/2$;

If $yx_{s+k}\notin E(G)$, then $od(y)\leq (p-k-1)/2$ and $id(y)\leq (p-k+1)/2$. In both cases we have $d(y)\leq p-k$, as required. (7) is proved.

From (4), (7) and condition ($M_0$) it follows that for every $l\in [1,k-1]$, 
$$ d(x_{s+l})\geq p+k-2.   \eqno (8)
   $$
Now for completes the proof of Claim 2, we will consider the cases $x_sy\in E(G)$ and $x_sy\notin E(G)$ separately.

\textbf{Case 1}. $x_sy\in E(G)$.

From (5) and (6) it follows that $yx_{s+k}\notin E(G)$ and $x_{s+k}y\in E(G)$. Thus we have
 $\{x_s,x_{s+k}\}\rightarrow y$. Notice that $Q:= x_{s+k}x_{s+k+1}\ldots x_{s}$ is a path of length $p-k-1$. Since (4) and (8), for every $l\in [1,k-1]$ we have 
$
 2(p-1)-|V(Q)|=p+k-2\leq d(x_{s+l})
$. 
Therefore, by Lemma 3.4, there exists a path $R:=y_1y_2\ldots y_{p-1}$ from $x_{s+k}$ to $x_s$ with vertex set $V(C_{p-1})$.
Hence, $D(p,2)=[y_1y; y_1y_2\ldots y_{p-1}y]$ is a Hamiltonian bypass, which contradicts our supposition that $G$ contains no  Hamiltonian bypass. 

\textbf{Case 2}. $x_sy\notin E(G)$.

From (5) it follows that $yx_s\in E(G)$. We may assume that $yx_{s+k}\notin E(G)$ (for otherwise in the converse digraph of $G$ we have $\{x_s,x_{s+k}\}\rightarrow y$, and hence in the converse digraph of $G$ the considered case $x_sy\in E(G)$ holds). This together with (5) implies that $x_{s+k}y\in E(G)$. Now using Lemma 3.5(i), 
Claim 1 and Case 1 ($x_sy, x_{s+k}y\in E(G)$), we may assume that $yx_{s+k+1}\in E(G)$ and $x_{s-1}y\in E(G)$. Then, by Lemma 3.5(iii), we have 
$$
x_{s+j}x_{s+j-1}\notin E(G)\quad \hbox{for all} \quad j\in [1,k]. \eqno (9)
$$
To be definite, assume that $x_1:=x_s$. Then $x_{s+k}=x_{k+1}$. Now we want to show that for any $i$ and $j$ with $1\leq i\leq j-1\leq k$,
$$
x_ix_j\in E(G) \quad \hbox{if and only if} \quad j=i+1. \eqno (10)
$$
Suppose, to the contrary, that this not so. Then for some  $i$ and $j$, $1\leq i\leq j-1\leq k$ and $j\not=i+1$, we have
$x_ix_j\in E(G)$. Consider the cycle $R:=x_ix_jx_{j+1}\ldots x_{k+1}yx_{k+2}\ldots x_{i-1}x_i$ of the length $p-j+i+1$.
By (8), for all $u\in \{x_{i+1},\ldots , x_{j-1}\}$ we have $ 2p-|V(R)|-1=p+j-i-2\leq p+k-2\leq d(u)$ since $j-i\leq k$. Therefore, by Lemma 3.3, there exists a cycle of length $p-1$ with vertex set $V(R)\cup \{x_{i+1},\ldots , x_{j-2}\}$ (if $j=i+2$, then $\{x_{i+1},\ldots , x_{j-2}\}=\emptyset$) that does not contain the vertex $x_{j-1}$. This contradicts Lemma 3.5(ii) since $d(x_{j-1})\geq p+k-2\geq p$ and $G$ contains no Hamiltonian bypass. Thus, (10) is true. 

From (9) and (10) it follows that for every subset $A\subseteq \{x_1,x_2,\ldots , x_{k+1}\}$ and for every 
$x_i\in \{x_1,x_2,\ldots , x_{k+1}\}$ the following holds
$$
d(x_i,A)\leq |A\setminus \{x_i\}|.   \eqno(11)
$$
In particular, from (11) it follows that
$$
d(x_1,\{x_2,x_3,\ldots , x_{k}\})\leq k-1 \quad \hbox{and} \quad d(x_{k+1},\{x_2,x_3,\ldots , x_{k}\})\leq k-1 .   \eqno(12)
$$
Put $P_1:=x_{k+2}x_{k+3}\ldots x_{p-1}x_{1}$ and $P_2:=x_{k+1}x_{k+2}\ldots x_{p-1}$ (possibly, $x_1=x_{k+2}$ and 
$x_{p-1}= x_{k+1}$). 
Now we want to show that the vertex $x_{k+1}$ cannot be inserted into $P_1$. Assume that this is not the case. Then there is a path,
say $Q$, from $x_{k+2}$ to $x_1$ with vertex set $\{x_{k+1},x_{k+2},\ldots , x_{p-1},x_{1}\}$. Then $|V(Q)|=p-k$. By (8) and (4), for every $u\in \{x_2,x_3,\ldots , x_{k}\}$ we have $d(u,V(G)\setminus \{y\})\geq p+k-2$. On the other hand, $2(p-1)-|V(Q)|=p+k-2$. 
Therefore, we can apply Lemma 3.4. In a result, we obtain an $(x_{k+2},x_{1})$-path, say $P$, with vertex set 
$V(G)\setminus \{y\}$. Then $D(p,2)=[yx_1;yP]$ is a Hamiltonian bypass since $yx_{k+2}\in E(G)$, which is a contradiction. This shows that
 $x_{k+1}$ cannot be inserted into $P_1$. 
Similarly, we can show that $x_1$ cannot be inserted into $P_2$.

Now using Lemma 3.2(ii) and the fact that $x_1x_{k+1}\notin E(G)$ (by (10)), we obtain
$$
d(x_{k+1},V(P_1))\leq p-k-1 \quad \hbox{and} \quad d(x_{1},V(P_2))\leq p-k-1 .   \eqno(13)
$$

 Assume first that $E(x_1,x_{k+1})=\emptyset$ (i.e., $x_1$ and $x_{k+1}$ are not adjacent). Then from condition $(M_0)$ it follows that  $d(x_1)+d(x_{k+1})\geq 2p-2$. 
This together with (12) and (13) implies that 
$
d(x_{1},V(P_2))= d(x_{k+1},V(P_1))= p-k-1
$. 
Therefore, since $E(x_1,x_{k+1})=\emptyset$, by Lemma 3.2(i) we have $x_1x_{k+2}\in E(G)$ and $x_{p-1}x_{k+1}\in E(G)$. Hence, $Q:=yx_1x_{k+2}\ldots x_{p-1}x_{k+1}y$ is a cycle of length $p-k+1$ and $2p-|V(Q)|-1=p+k-2$. On the other hand, by (8), for every $u\in \{x_2,x_3,\ldots , x_k\}$ the following holds $d(u)\geq 2p-|V(Q)|-1$.
Therefore, by Lemma 3.3, there exists a cycle of length $p-1$ with vertex set $V(G)\setminus \{x_k\}$. This contradicts Lemma 3.5(ii) since $d(x_{k})\geq p+k-2\geq p$ and $G$ has no Hamiltonian bypass.

 Assume now that $E(x_1,x_{k+1})\not=\emptyset$. Then from (10) we have $x_1x_{k+1}\notin E(G)$ and hence, 
$x_{k+1}x_1\in E(G)$. Using (8), (11) and the fact that $p\geq 3$, it is easy to see that $k\not= p-2$. From (4), (8) and (10) it follows  
$
d(u,V(G)\setminus \{x_{k+1},y\})\geq p+k-3 \quad  \hbox{for all} \quad u\in \{x_2,x_3,\ldots , x_k\}
$. 
Hence,
 $$
 2(p-2)-|V(P_1)|=p+k-3\leq d(u,V(G)\setminus \{x_{k+1},y\},) \quad  \hbox{for all} \quad u\in \{x_2,x_3,\ldots , x_k\}.
$$ 
 Now we can apply Lemma 3.4 to the path $P_1$ and to the set $\{x_2,x_3,\ldots , x_k\}$.
In a result, we obtain an $(x_{k+2},x_1)$-path,  say $H$, with vertex set 
$V(C)\setminus \{x_{k+1},y\}$. Then $D(p,2)=[x_{k+1}x_1;x_{k+1}yH]$ is a Hamiltonian bypass  since the arcs $x_{k+1}y$,    $yx_{k+2}$ are in $E(G)$, which is a contradiction. This contradiction completes the proof of Claim 2. \fbox \\\\

Now we are ready to finish the proof of Lemma 4.2.
From Claim 2 and Lemma 3.5(ii) it follows that $p$ is odd and $id(y)=od(y)=(p-1)/2$. Using Lemma 3.5(i), we may assume that 
$O(y)=\{x_2,x_4,\ldots , x_{p-1}\}$. It is easy to see that $p\geq 5$  since $G$ is not isomorphic to  $D_{p-1,1}$. Lemma 3.5(i) and Claim 1 imply that 
 $I(y)=\{x_1,x_3,\ldots , x_{p-2}\}$. Therefore, by Lemma 3.5(iii), we have  that $x_{i+1}x_i\notin E(G)$ for all $i\in [1,p-1]$ since $G$ contains no Hamiltonian bypass.  

 We first consider  the case when for some $i\in [1,p-1]$, $x_{i-1}x_{i+1}\in E(G)$. Without loss of generality, we may assume that $i$ odd ,i.e., $x_iy, yx_{i-1}$ and $yx_{i+1}\in E(G)$.
(for otherwise, we will consider the converse digraph of $G$). If $x_{i}x_{i+2}\in E(G)$, then 
$D(p,2)=[x_{i-1}x_{i+1};x_{i-1}x_ix_{i+2}\ldots x_{i-2}yx_{i+1}]$, which is a contradiction. We may therefore assume that
$x_{i}x_{i+2}\notin E(G)$. 
If $x_{i-2}x_{i}\in E(G)$ and $p\geq 7$, then the cycle 
$x_{i-1}x_{i+1}x_{i+2}yx_{i+3}\ldots \\ x_{i-2}x_{i-1}$ has length $p-1$ and does not contain $x_i$, but  $|E(\{x_{i-2},x_{i-1}\}\rightarrow x_{i})|=2$,
which contradicts Lemma 3.5(i). If $x_{i-2}x_{i}\in E(G)$ and $p=5$, then it is easy to check that $x_{i+1}x_{i-1}\notin E(G)$ and
$x_{i}x_{i-2}\notin E(G)$. Therefore, $G$ is isomorphic to $ T_5$, which contradicts our supposition. (To see this, we assume that $x_i=x_1$, then $x_{i-1}= x_4$, $x_{i+1}=x_2$, $x_{i-2}=x_3$, and consider the following mapping: $x_1\mapsto z_3$, $x_2\mapsto z_4$, $x_3\mapsto z_1$, $x_4\mapsto z_2$). Thus we have proved that
$$
\hbox{if}\quad x_{i-1}x_{i+1}\in E(G),\quad  \hbox{then} \quad x_{i}x_{i+2}\notin E(G) \quad \hbox{and} \quad x_{i-2}x_{i}\notin E(G). \eqno (14)
$$
If there exists an $(x_{i+2},x_{i-2})$-path, say $Q$, with vertex set $V(C_{p-1})\setminus \{x_{i+1},x_{i-1}\}$, then 
$D(p,2)=[yx_{i-1}; yx_{i+1}Qx_{i-1}]$, which is a contradiction. We may  therefore assume that there is no 
 $(x_{i+2},x_{i-2})$-path  with vertex set $V(C_{p-1})\setminus \{x_{i+1},x_{i-1}\}$. This means that the vertex $x_i$ cannot 
be inserted into  $x_{i+2}x_{i+3}\ldots x_{i-2}$. Now using Lemma 3.2, (14) and the fact that
$|E(x_i,\{x_{i-1},x_{i+1},y\}|=3$, we obtain that $d(x_{i})\leq p-2$, which contradicts Claim 2 since
$x_{i-2}yx_{i-1}x_{i+1}\ldots x_{i-2}$ is a cycle of length $p-1$, which does not contain $x_i$, but $d(x_{i})\leq p-2$.

 We next consider the case when for all $i\in [1,p-1]$,  
$x_{i-1}x_{i+1}\notin E(G)$.  
It is easy to check that $p\geq 7$. Indeed, if $p=5$, then $E(x_2,x_4)=\emptyset$ and $d(x_2)=d(x_4)=3$, i.e.,  
$d(x_2)+d(x_4)=6$, which contradicts condition $(M_0)$. Thus, $p\geq 7$.

Assume that $E(x_{i-1},x_{i+1})=\emptyset$ for some $i\in [1,p-1]$. Then, by condition $(M_0)$, 
$$
d(x_{i-1})+d(x_{i+1})\geq 2p-2.    \eqno (15)
$$ 
To be definite, assume that $i$ even. We claim that the path $x_{i+2}x_{i+3}\ldots x_{i-1}$ 
 cannot be extended with the vertex $x_{i+1}$. For otherwise there is an $(x_{i+2},x_{i-1})$-path, say $R$, with vertex set $\{x_{i+1},x_{i+2}, \ldots , x_{i-1}\}$. Then, $yRx_{i-1}y$ is a cycle of length $p-1$, which does not contain $x_i$, but $\{x_{i-1},y\}\rightarrow x_i$, which contradicts Lemma 3.5(i).  Similarly, the path $x_{i+1}x_{i+2}\ldots x_{i-2}$  cannot be extended with the vertex $x_{i-1}$.
 Now using Lemma 3.2 and the facts that $d(x_{i+1},\{y,x_{i},x_{i+2}\})=3$, $E(x_{i-1},x_{i+1})=\emptyset$ and $x_{i+1}x_{i+3}\notin E(G)$ (by our assumption), we obtain 
$$
d(x_{i+1})=d(x_{i+1},\{y,x_{i},x_{i+2}\})+ d(x_{i+1},\{x_{i+3},x_{i+4},\ldots , x_{i-2}\})\leq p-2.
$$
  Similarly, $d(x_{i-1})\leq p-2$. 
The last two inequalities contradict (15).

Assume now that for all $i\in [1,p-1]$, $E(x_{i-1},x_{i+1})\not=\emptyset$. Then $x_{i+1}x_{i-1}\in E(G)$ for all $i\in [1,p-1]$.
Hence,  $D(p,2)=[x_{4}x_{2};x_{4}x_{5}x_{3}yx_6\ldots x_{p-1}x_{1}x_{2}]$, a contradiction. This contradicts our supposition, and completes the discussion of Case 2.
 Lemma 4.2 is proved. \fbox \\\\

\textbf{Proof of Lemma 4.3}.  

 Let $G$ be a digraph of order $p\geq 3$ satisfying the conditions of Lemma 4.3. By contradiction, suppose that   $G$ is not isomorphic to  $C_3$ and contains no Hamiltonian bypass.  Let 
$C_p:=x_1x_2\ldots x_{p}x_1$ be an arbitrary Hamiltonian cycle in $G$. It is easy to see that $x_{i-1}x_{i+1}\notin E(G)$
(since $G$ contains no cycle of length $p-1$) and $x_ix_{i-1}\notin E(G)$ for all $i\in [1,p]$.  
Using these, it is not difficult to show  that $p\geq 6$. We first prove the following claim.

\textbf{Claim 3}. For all $i\in [1,p]$, the vertices  $x_{i-1}$ and $x_{i+1}$ are not adjacent.

{\it Proof}. By contradiction, assume that  $x_{i-1}$ and $x_{i+1}$ are  adjacent  for some $i\in [1,p]$.
Since $x_{i-1}x_{i+1}\notin E(G)$, we have  $x_{i+1}x_{i-1}\in E(G)$. To be definite, assume that $x_i=x_1$, i.e.,
$x_2x_p\in E(G)$. We distinguish two cases.

\textbf{Case 1.} There exists an integer $k\in [4,p-1]$ such that $x_1x_k\in E(G)$.

Assume that $k$ is the minimum with this property, i.e.,
$$
E(x_1\rightarrow \{x_3,x_4,\ldots , x_{k-1}\})=\emptyset. \eqno (16)
$$
Assume first that there exists an integer $s\in [3,k-1]$ such that $x_sx_1\in E(G)$. Since $G$ contains no Hamiltonian bypass, it follows that $s\leq k-2$.  Assume that  $s$ is maximal with this property. Then by (16), we have
$$
E(x_1, \{x_{s+1},x_{s+2},\ldots , x_{k-1}\})=\emptyset. \eqno (17)
$$
Put $P_1:=x_kx_{k+1}\ldots x_{p-1}$; $P_2:=x_3x_{4}\ldots x_{s}$; $P_3:=x_2x_{3}\ldots x_{s}x_1x_kx_{k+1}\ldots x_p$ and 
$A:=\{x_{s+1},x_{s+2},\ldots ,$ $ x_{k-1}\}$; $a:=|A|=k-s-1$. Since the paths $P_1$ and  $P_2$ cannot be extended with the vertex $x_1$ and $x_1x_3\notin E(G)$, $x_{p-1}x_1\notin E(G)$, from Lemma 3.2(ii) and (17) it follows that
$$
d(x_1)= d(x_1,V(P_1))+d(x_1,V(P_2))+d(x_1,\{x_2,x_p\})\leq |V(P_1)|+|V(P_2)|+2=p-a-1.
$$
From this, (17) and condition ($M_0$) we get that for all $u\in A$,
$$
d(u)\geq p+a-1. \eqno (18)
$$
On the other hand, it is clear that the path $P_3$ cannot be extended with all the vertices of $A$ since $x_2x_p\in E(G)$. Therefore, for some vertices 
$u_1,u_2,\ldots , u_d$ of $A$ , where $1\leq d\leq a$, by Lemma 3.2, the following holds $d(u_i)\leq p+d-1\leq p+a-1$.
This together with (18) implies that $d=a$ and the induced subdigraph $G\langle A\rangle $ is a complete digraph. Therefore, $a=1$, (i.e., $s=k-2$)
since $x_ix_{i-1}\notin E(G)$ for all $i\in [1,p]$. It is not difficult to see that the vertex $x_{k-1}$ cannot be inserted
neither into   $x_2x_3\ldots x_{k-2}$ nor  $x_kx_{k+1}\ldots x_p$. This together with
$$
E(x_{k-1}\rightarrow \{x_{k+1},x_{k-2}\})=E(\{x_{k-3},x_{k}\}\rightarrow x_{k-1})=\emptyset 
$$
and Lemma 3.2(ii) implies that $d(x_{k-1})\leq p-2$, which contradicts (18).

Assume next that $E(\{x_3,x_4, \ldots , x_{k-1}\}\rightarrow x_{1})=\emptyset$. This together with (16) gives
$$
E(x_1,\{x_3,x_4, \ldots , x_{k-1}\})=\emptyset.   \eqno (19)
$$
Therefore, since the path $P_1=x_kx_{k+1}\ldots x_{p-1}$ cannot be extended with $x_1$ and $x_{p-1}x_1\notin E(G)$ using Lemma 3.2(ii), we obtain  
$
d(x_1)=d(x_1,V(P_1))+d(x_1,\{x_2,x_p\})\leq p-k+2
$
(note that $|V(P_1)|=p-k$). Now from (19) and condition ($M_0$) it follows that for every $i\in [3,k-1]$, $d(x_i)\geq p+k-4$. Now we will
consider the cycle 
$C_{p-k+2}:=x_1x_kx_{k+1}\ldots x_px_1$ of length $p-k+2$. It is clear that $C_{p-k+2}$ cannot be extended with all the vertices 
$x_3,x_4,\ldots , x_{k-1}$ (for otherwise $G$ contains a cycle of length $p-1$).
 Therefore, for some $u_1,u_2,\ldots , u_d \in \{x_3,x_4,\ldots , x_{k-1}\}$, where $1\leq d\leq k-3$, by Lemma 3.1  the following holds
$$
d(u_i)=d(u_i,V(G)\setminus \{u_1,u_2,\ldots , u_d,x_2\})+d(u_i, \{u_1,u_2,\ldots , u_d,x_2\})\leq p+d-1.
$$ 
Therefore, $p+k-4\leq d(u_i)\leq p+d-1$. This implies that $d=k-3$ and the induced  subdigraph $G\langle \{x_2,x_3, \ldots , x_{k-1}\}\rangle$ is a complete digraph. Then $x_3x_2\in E(G)$,  which is a contradiction.

\textbf{Case 2.} For all $i\in [3,p-1]$, $x_1x_i\notin E(D)$.

We may assume that  $E(\{x_3,x_4,\ldots , x_{p-2}\}\rightarrow x_1)=\emptyset$ (for otherwise in the converse digraph of $G$ 
we have the considered Case 1). Therefore
$$
E(x_1, \{x_3,x_4,\ldots , x_{p-1}\})=\emptyset. \eqno (20)
$$
This together with $x_2x_1\notin E(G)$ and $x_1x_p\notin E(G)$ implies that $d(x_1)=2$. Hence, by (20) and condition 
($M_0$) for all $i\in [3,p-1]$ we have $d(x_i)\geq 2p-4$. On the other hand, since (20) and $x_{i+1}x_i\notin E(G)$, we have $d(x_i)\leq 2p-6$, where $i\in [3,p-1]$. Thus we have a contradiction.  
This  completes the proof of Claim 3.  \fbox \\\\

Using Claim 3, condition ($M_0$) and the fact that $x_ix_{i-1}\notin E(G)$, it is not difficult to show that $p\geq 8$. It is clear that $G$ is not a directed cycle. We choose an arc $x_jx_k\in E(G)$ with  $k\not=j+1$ such that $|\{x_j,x_{j+1},\ldots , x_k\}|$ is the smallest possible. To be definite, assume
that $j=1$. By Claim 3, $4\leq k\leq p-2$. We claim  that for all $i\in [k+1,p]$ the following holds
$$
|E(x_2\rightarrow x_i)|+|E(x_{k-1}\rightarrow x_{i+1})|\leq 1 \quad \hbox{and} \quad |E(x_{i}\rightarrow x_{k-1})|+
|E(x_{i-1}\rightarrow x_{2})|\leq 1.   \eqno (21) 
$$
Indeed, if $x_2x_i\in E(G)$ and $x_{k-1}x_{i+1}\in E(G)$, then
$D(p,2)=[x_2x_i;x_2x_3\ldots x_{k-1}x_{i+1}\ldots x_p$ $ x_1x_k\ldots x_i]$;
if $x_{i-1}x_2\in E(G)$ and $x_{i}x_{k-1}\in E(G)$, then
$D(p,2)=[x_ix_{k-1};x_ix_{i+1}\ldots x_{p}x_{1}x_k$ $\ldots x_{i-1}x_2\ldots x_{k-1}]$, which contradicts the supposition that $G$ contains no Hamiltonian bypass. 

From the first inequality of (21), Claim 3, the minimality of $k$ and the fact that for all $i\in [1,p]$, $x_ix_{i-1}\notin E(G)$  it follows that
$$
od(x_2)+od(x_{k-1})=od(x_2,\{x_p,x_1,x_2,\ldots , x_k\})+ od(x_{k-1},\{x_1,x_2,\ldots , x_{k-1}, x_k, x_{k+1}\})
$$ 
$$
+od(x_2,\{x_{k+1},x_{k+2},\ldots , x_{p-1}\})+od(x_{k-1},\{x_{k+2},x_{k+3},\ldots , x_p\})$$ $$\leq 1+k-3 
 + \sum^{p-1}_{j=k+1}(|E(x_2\rightarrow x_j)|+|E(x_{k-1}\rightarrow x_{j+1})|)\leq p-3.
$$
Similarly, using the second inequality of (21), we obtain
$$
id(x_2)+id(x_{k-1})=id(x_2,\{x_p,x_1,x_2,\ldots ,x_{k-1}, x_k\})+ id(x_{k-1},\{x_1,x_2,\ldots , x_{k-1}, x_k, x_{k+1}\})
$$ 
$$
+ id(x_2,\{x_{k+1},x_{k+2},\ldots , x_{p-1}\})+id(x_{k-1},\{x_{k+2},x_{k+3},\ldots , x_p\})
 $$
$$
 \leq k-3+1+ \sum^{p-1}_{j=k+1}(|E(x_j\rightarrow x_2)|+|E(x_{j+1}\rightarrow x_{k-1})|)\leq k-2+p-k-1=p-3.
$$
Therfore,
$
d(x_2)+d(x_{k-1})\leq 2p-6. 
$
Now, taking into account condition $(M_0)$, we obtain  that $x_2$ and $x_{k-1}$ are adjacent and $d(x_2)\leq p-3$ or $d(x_{k-1})\leq p-3$. We may assume that $d(x_2)\leq p-3$ (for otherwise we consider the converse digraph of $G$). We will consider the cases $k\geq 5$ and $k=4$ separately. 

\textbf{Case 1.} $k\geq 5$.

Since $x_2$ and $x_{k-1}$ are adjacent, from the minimality of $k$ it follows that $x_{k-1}x_2\in E(G)$. This and Claim 3 imply that $k\geq 6$. Since  $x_{k-1}x_2\in E(G)$ and the vertices $x_2$,  $x_{4}$ are nonadjacent (Claim 3), from the minimality of $k$ it follows that there exists an integer $s\in [5,k-1]$ such that $x_sx_2\in E(G)$ and 
$
E(\{x_3,x_4,\ldots , x_{s-1}\}\rightarrow x_2)=\emptyset
$.

Put $A:=\{x_3,x_4,\ldots , x_{s-1}\}$. From the minimality of $k$ we have
$$
E(\{x_1,x_p\}\rightarrow A)=\emptyset \quad \hbox{and} \quad E(A\rightarrow \{x_{s+1},x_{s+2},\ldots , x_{k}\}
)= \emptyset \eqno (22)
$$
and for any subset $A_1\subseteq A$ and  for every vertex $u\in A_1$ the following holds
$$
d(u,\{x_s,x_{s+1},\ldots , x_{k-1}\})\leq |\{x_s,x_{s+1},\ldots , x_{k-1}\}|; \eqno (23)
$$
$$
d(u,A_1)\leq |A_1|-1, \quad d(u,\{x_2\})\leq 1 \quad \hbox{and} \quad d(u,\{x_k\})\leq 1. \eqno (24)
$$
We extend the path $x_kx_{k+1}\ldots x_p$ with the vertices of $A$ as much as possible. It is clear that some vertices 
$u_1,u_2,\ldots , u_d\in A$, $1\leq d\leq |A|$, do not on the obtained extended path, say $R$. Notice that
$|R|=p-k+s-d-2$. Using (22)-(24) and Lemma 3.2(iii), we obtain
$$
d(u_i)=d(u_i,R)+d(u_i,\{x_s,x_{s+1},\ldots , x_{k-1}\})+d(u_i,\{x_1,x_2\})+d(u_i,\{u_1,u_2,\ldots ,u_d\})
$$
 $$
\leq |R|-1+k-s+2+d-1=p-2.
$$
This together with $d(x_2)\leq p-3$ and condition ($M_0$) implies that the vertex $x_2$ and  every vertex $u_i$    are adjacent. Therefore from the minimality of $s$ and $k$ it follows that $d=1$, $u_1=x_3$ and there is a path
 $Q:=y_1y_2\ldots y_{p-1}$ from $x_s$ to $x_2$ with vertex set  $V(G)\setminus \{x_3\}$. Notice that
 $y_1=x_s$, $y_2=x_{s+1}$, $x_1=y_{p-2}$, $x_2=y_{p-1}$, $x_p=y_{p-3}$, $x_k=y_{k-s+1}$ and $y_1y_{p-1}\in E(G)$. It is clear that $E(x_3\rightarrow V(Q))\not=\emptyset$ since $x_3 x_4\in E(G)$. Therefore, for some $l\in [3,p-3]$, $x_3y_l\in E(G)$ since $E(x_3\rightarrow \{x_1,x_2,x_s,x_{s+1}\})=\emptyset$. Let $l$ be the smallest with these properties, i.e.,
$$
E(x_3\rightarrow \{y_1,y_2,\ldots , y_{l-1}\})=\emptyset. \eqno (25)
$$

Assume first that there exists a $q\in [1,l-1]$ such that $y_qx_3\in E(G)$. Notice that $q\leq l-2$ since $y_1y_{p-1}\in E(G)$. Let $q$ be the maximum with this property,  i.e., 
$
E(\{y_{q+1},\ldots , y_{l-1}\}\rightarrow x_3)=\emptyset 
$.
This together with (25) implies that 
$$
E(x_3,\{y_{q+1},y_{q+2},\ldots , y_{l-1}\})=\emptyset. \eqno (26)
$$
Now put $P_1:=y_1y_2\ldots y_q$ and $P_2:=y_ly_{l+1}\ldots y_{p-3}$. Since the paths $P_1$ and $P_2$ cannot be extended with 
 $x_3$ and $x_3y_1\notin E(G)$, $y_{p-3}x_3\notin E(G)$ (by the minimality of $k$) and $E(x_3,x_{p-2})=\emptyset$, using Lemma 3.2(ii) and (26), we obtain
$$
d(x_3)=d(x_3,V(P_1))+d(x_3,V(P_2))+d(x_3,\{y
_{p-1}\})
$$
$$
\leq |V(P_1)|+|V(P_2)|+1=p+q-l-1. \eqno (27)
$$
We also have that the path $R:=y_1y_2\ldots y_qx_3y_ly_{l+1}\ldots y_{p-1}$ cannot be extended with all the vertices 
$y_{q+1},$ $y_{q+2},\ldots , y_{l-1}$. Therefore, by Lemma 3.2, for some vertices 
$u_1,u_2,\ldots , u_d\in \{y_{q+1},y_{q+2},\ldots , y_{l-1}\}$, where  $1\leq d\leq l-q-1$, the following holds
$$
d(u_i)\leq p+d-1\leq p+l-q-2. \eqno (28)
$$
Combining this together with (27), we obtain that $d(x_3)+d(u_i)\leq 2p-3$, which contradicts condition ($M_0$) since the vertices $x_3$ and $u_i$ are not adjacent.

Assume next that $E(\{y_1,y_2,\ldots , y_{l-1}\}\rightarrow x_3)=\emptyset$. Then $E(x_3,\{y_1,y_2,\ldots , y_{l-1}\})=\emptyset$ since $l$ is minimal. Therefore, 
$$
d(x_3)=d(x_3,V(P_2))+d(x_3,\{y_{p-1}\})\leq |V(P_2)|+1=p-l-1. 
$$
Then, by condition ($M_0$),  for all $j\in [1,l-1]$, $d(y_j)\geq p+l-1$.  Now consider the cycle 
$R:=x_3y_ly_{l+1}\ldots y_{p-1}x_3$ (recall that $y_{p-1}=x_2$) of length $p-l+1$. It is easy to check that
 $2p-|V(R)|-1=p+l-2\leq d(y_i)$. Therefore, by Lemma 3.3, there exists a cycle of length $p-1$, 
which contradicts the assumption of Lemma 4.3 that $G$ contains no cycle of length $p-1$. This completes the discussion of case $k\geq 5$. 

\textbf{Case 2.} $k=4$, i.e., $x_1x_4\in E(G)$.

Since $G$ contains no cycle of length $p-1$, it follows that $x_2$ ($x_3$) cannot be inserted into the cycle 
$C_{p-2}:= x_1x_4x_5\ldots x_px_1$. Recall that $E(x_2,\{x_4,x_p\})=\emptyset$ (Claim 3) and $d(x_2)\leq p-3$ (by our assumption). 
This together with condition $(M_0)$ implies that $d(x_p)\geq p+1$ and $d(x_4)\geq p+1$.
It is not difficult to show that $x_3x_6\notin E(G)$. Indeed, if $x_3x_6\in E(G)$, then $x_3x_6x_7\ldots x_px_1x_2x_3$ is a cycle of length $p-2$. Since $d(x_4)\geq p+1$ and $x_5x_4\notin E(G)$, from Lemma 3.1 it follows that $x_4$ can be inserted into this cycle, i.e., $G$ contains a cycle with vertex set $V(G)\setminus \{x_5\}$, which contradicts the assumption of Lemma 4.3 that $G$ contains no cycle of length $p-1$.

 \textbf{Subcase 2.1}. For some $l\in [7,p]$, $x_3x_l\in E(G)$ (recall that $p\geq 8$).

Let $l$ be the minimum with this property, i.e.,
$
E(x_3\rightarrow \{x_5,x_6,\ldots , x_{l-1}\})=\emptyset
$.
Let $E(\{x_6,x_7,\ldots , x_{l-1}\}\rightarrow x_3)\not=\emptyset$. Observe that $x_{l-1}x_3\notin E(G)$, since $G$ contains no cycle of length $p-1$. Then there exists an integer $s\in [6,l-2]$ such that 
$x_sx_3\in E(G)$ and
$$
E(x_3, \{x_{s+1},x_{s+2},\ldots , x_{l-1}\})=\emptyset. \eqno (29)
$$

Put $P_1:=x_6x_7\ldots x_s$ and $P_2:=x_lx_{l+1}\ldots x_{p}$. Since $G$ contains no cycle of length $p-1$, it follows that $x_3$ cannot be inserted neither into $P_1$ nor in $P_2$. Now using Lemma 3.2, (29) and the facts that 
$E(x_3,\{x_5,x_1\})=\emptyset$ and $x_3x_6\notin E(G)$, we obtain
$$
d(x_3)=d(x_3,V(P_1))+d(x_3,V(P_2))+d(x_3,\{x_2,x_4\})\leq p-l+s-1. 
$$
Therefore, because of (29) and condition ($M_0$), for every $u\in \{x_{s+1},x_{s+2},\ldots , x_{l-1}\}$ we have
$$
d(u)\geq p+l-s-1. \eqno (30)
$$
Since $G$ contains no cycle of length $p-1$, it follows that the cycle 
$x_1x_4x_5\ldots x_sx_3x_lx_{l+1}\ldots$ $ x_px_1$
cannot be extended with all the vertices of $\{x_{s+1},x_{s+2},\ldots , x_{l-1}\}$. This means that for some vertices 
$u_1,u_2,\ldots , u_d\in \{x_{s+1},x_{s+2},\ldots , x_{l-1}\}$, $1\leq d\leq l-s-1$, by Lemma 3.1, the following holds
$$
d(u_i)=d(u_i,V(G)\setminus \{x_2,u_1,u_2,\ldots ,u_d\})+d(u_i,\{x_2,u_1,u_2,\ldots ,u_d\})\leq p+d-1\leq p+l-s-2,
$$
 which contradicts (30).

Let now $E(\{x_6,x_7,\ldots , x_{l-1}\}\rightarrow x_3)=\emptyset$. Then, by  minimality of $l$,   
$
E(x_3, \{x_{5},x_{6},\ldots ,\\ x_{l-1}\})=\emptyset
$,
 and by Lemma 3.2,  $d(x_3)=d(x_3, V(P_2))+d(x_3,\{x_2,x_4\})\leq p-l+4$. Since the vertices $x_3$ and $x_i$, $i\in [5, l-1]$ are not adjacent, the last equality and the last inequality
together with condition $(M_0)$ imply that for every $i\in [5,l-1]$,
$d(x_i)\geq p+l-6
$. 
Consider the cycle $Q:=x_3x_lx_{l+1}\ldots x_px_1x_2x_3$ of length $p-l+4$. It is easy to see that the cycle $Q$ cannot be extended with all the vertices of  $\{x_5,x_6,\ldots , x_{l-1}\}$ since $G$ has no cycle of length $p-1$. This means that for some vertices $u_1,u_2,\ldots u_d$ of  $\{x_5,x_6,\ldots , x_{l-1}\}$, where $1\leq d\leq l-5$, by Lemma 3.1, the following holds 
$$
d(u_i)=d(u_i,V(G)\setminus \{x_4,u_1,u_2,\ldots , u_d\})+d(u_i,\{x_4,u_1,u_2,\ldots , u_d\})\leq p+d-1\leq p+l-6.
$$
From this and $d(x_i)\geq p+l-6$ it follows that $d=l-5$ and $d(u_i, \{x_4\})=2$, in particular, $x_5x_4\in E(G)$, which is a contradiction since $G$ has no Hamiltonian bypass.

 \textbf{Subcase 2.2}. $E(x_3\rightarrow \{x_6,x_7,\ldots , x_{p}\})=\emptyset$.

Let for some $l\in [5,p]$, $x_lx_3\in E(G)$. Pick a maximum such $l$. Then $l\geq 6$ (Claim 3). We have that 
 $
E(x_3, \{x_{l+1},x_{l+2},\ldots , x_{p},x_1\})=\emptyset
$.
Since $x_3$ cannot be inserted into the path $x_6x_7\ldots x_l$ and $x_3x_6\notin E(G)$, using Lemma 3.2(ii), we obtain 
$$
d(x_3)=d(x_3,\{x_6,x_7,\ldots ,x_l\})+d(x_3,\{x_2,x_4\})\leq l-3.
$$
 This together with condition ($M_0$) implies that for every $x_i\in \{x_1,x_{l+1},x_{l+2},\ldots , x_p\}$ the following holds
$d(x_i)\geq 2p-l+1$ and $d(x_i,V(G)\setminus \{x_2\})\geq 2p-l-1$. 
Now we consider the cycle $Q:=x_3x_4\ldots x_lx_3$ of length $l-2$. It is easy to see that 
$$
2(p-1)-|V(Q)|-1= 2p-l-1\leq d(x_i,V(G)\setminus \{x_2\}).
$$  
Therefore, we can apply Lemma 3.3 to the cycle $Q$. In a result, we obtain a cycle of length $p-1$ with vertex set $V(G)\setminus \{x_2\}$, which is a contradiction.

Assume finally that $E(x_3,\{x_6,x_7,\ldots , x_{p}\})=\emptyset$. Then $d(x_3)=2$ since $d(x_3,\{x_1,x_5\})=0$,
$x_3x_2\notin E(G)$ and  
$x_4x_3\notin E(G)$. This together with condition ($M_0$) implies that $d(x_1)\geq 2p-4$.  On the other hand, since $d(x_1,\{x_3, x_{p-1}\})=0$, $x_2x_1\notin E(G)$ and $x_1x_p\notin E(G)$, we have that $d(x_1)\leq 2p-8$, a contradiction.   This contradiction completes the discussion of case $k=4$. Lemma 4.3 is proved. \fbox \\\\

\textbf{Proof of Lemma 4.4}.

 Let $G$ be a digraph of order $p\geq 3$ satisfying the conditions of Lemma 4.4. Then, $p\geq 4$. Suppose that $G$ is not isomorphic to $D_{p-k,k}$, where $k\in [1, p-2]$. Let $C_m:=x_1x_2\ldots x_mx_1$ be a longest cycle
in $G$. Then, $2\leq m\leq p-2$. Let $D_1$, $D_2$ ,..., $D_s$ be the strong components of 
$G\langle V(G)\setminus V(C_m)\rangle$ labelled in such a way that no vertex of $D_i$ dominates a vertex of $D_j$ whenever
$i>j$. By Theorem 1.9,  

(i) for every $k\in [1,s]$, $D_k$ is a complete digraph; 

(ii) $V(G)\setminus V(C_m)$ contains a Hamiltonian path, say $y_1y_2\ldots y_{p-m}$ ;

 (iii) there exist two distinct vertices $x_k$ and $x_j\in V(C_m)$ (to be definite, we assume that $x_j=x_1$) such that $k\leq m-1$, $x_k\rightarrow V(G_1)\rightarrow x_1$, 
$
E(\{x_{k+1},x_{k+2},\ldots , x_m\},V(D_1))=\emptyset 
$ and $d(u,V(C_m))=k+1$ for all $u\in V(G_1))$. Moreover, if $s\geq 2$, then $\{x_{k+1},x_{k+2},\ldots ,$ $ x_m\}\rightarrow V(D_s)$. Therefore, if 
$s\geq 2$, then   $D(p,2)=[x_{k+1}y_{p-m}; x_{k+1}x_{k+2}\ldots x_ky_1y_2\ldots y_{p-m}]$ is a Hamiltonian bypass. We may therefore assume  that $s=1$. Then, $|V(D_1)|\geq 2$. This together with $d(u,V(C_m))=k+1$  implies that there exists an integer $l$,  $1\leq l\leq k$, such that 
$
\{x_{l},x_{l+1}\ldots , x_k\}\rightarrow V(D_1)\rightarrow \{x_{1},x_2, \ldots , x_{l}\}  
$. From the above observations we have: if   $l\leq k-1$, then $D(p,2)=[x_{k}y_{p-m}; x_{k}x_{k+1}\ldots x_mx_1\ldots x_{k-1}y_1y_2\ldots y_{p-m}]$,  and if
$l=k$, then  $D(p,2)=[y_1x_1; y_1y_2\ldots y_{p-m}x_2x_3\ldots x_mx_1]$. Thus in all cases we have shown that $D$ contains a 
Hamiltonian bypass. Lemma 4.4 is proved, and completes the proof of Theorem 1.18. \fbox \\\\

From Theorem 1.18 it follows the following corollaries 1 and 2.
 
\textbf{Corollary 1} (\cite{[1]}). {\it Let $G$ be a digraph of order $p\geq 3$. If $od(x)+id(y)\geq p$ for all pairs of vertices $x$ and $y$ 
such that there is no arc from $x$ to $y$, then $G$ contains a Hamiltonian bypass.}

\textbf{Corollary 2} (\cite{[1]}). {\it Let $G$ be a digraph of order $p\geq 3$ with minimum degree at least $p$.  
Then $G$ contains a Hamiltonian bypass.}

\textbf{Corollary 3}. {\it Let $G$ be a digraph of order $p\geq 3$ satisfying condition $(M_1)$.   
Then $G$ contains a Hamiltonian bypass unless $G\in \{C_3, T_5\}$.}

\textbf{Proof}. It is easy to check that if a digraph $G$ sutisfies condition $(M_1)$, then it is
  neither isomorphic to a digraph of type $D_0$ nor to a digraph of type $D_{p-k,k}$. For $p=3$, the theorem clearly is true. Assume that $p\geq 4$. In order to prove the corollary, by Theorem 1.18, it suffices to consider the case when $G$ is not strongly connected. Let $G_1, G_2,\ldots , G_s$, $s\geq 2$, be the strongly connected components of $G$ labelled in such a way that no vertex of $G_i$ dominates a vertex of $G_j$ whenever $i>j$. It is not difficult to show that for any pair of integers $i, j$, $1\leq i<j\leq s$, there exist vertices 
$x\in V(G_i)$ and $y\in V(G_j)$ such that $xy\in E(G)$. Indeed, in the converse case for any  $u\in V(G_i)$ and $v\in V(G_j)$ we have 
$$
d(u)+d(v)=d(u,V(G_i))+d(u,V(G)\setminus \{V(G_i)\cup V(G_j)\})+d(v,V(G_j))$$ $$+d(v,V(G)\setminus \{V(G_i)\cup V(G_j)\})$$ $$ \leq 
2|V(G_i)|-2+2p-2|V(G_i)|-2|V(G_j)|+2|V(G_j)|-2=2p-4,
$$
which contradicts condition $(M_1)$. In particular, there are two vertices $x\in V(G_1)$ and $y\in V(G_s)$ such that
 $xy\in E(G)$. If in $G$  instead of the $xy$ we replace the arc $yx$, then we obtain a strong digraph, say $G'$. By the Meyniel theorem, $G'$ contains a Hamiltonian cycle containing the arc $yx$, which in turn implies that $G$ contains a Hamiltonian bypass. The corollary is proved. \fbox \\\\

\textbf{Proof of Theorem 1.19}. 

Let $G$ be a strong digraph of order $p\geq 4$ satisfying condition $(M_1)$.

Assume first that $G$ contains a cycle of length $p-1$. Let $C_{p-1}:=x_1x_2\ldots x_{p-1}x_1$ be a cycle of length $p-1$ in $G$ and let $y$ be the vertex that is not on $C_{p-1}$. If $yx_i$ and $x_{i+1}y\in E(G)$, then $D(p,3)=[x_{i+1}yx_i; x_{i+1}x_{i+2}\ldots x_i]$. We may therefore assume that 
$$
|E(y\rightarrow x_i)|+|E(x_{i+1}\rightarrow y)|\leq 1 \quad \hbox{for all} \quad i\in [1,p-1]. \eqno (31)
$$
Since  $G$ is strong, from condition ($M_1$) it follows that there exist two distinct vertices $x_i$ and $x_j$
such that $yx_i$ and $x_jy\in E(G)$. Therefore, we can choose two  distinct integers $k, l\in [1,p-1]$  such that $x_{k}y\in E(G)$ (we may assume that $x_k=x_1$),  
$yx_l\in E(G)$ and 
$$
E(y,\{x_{l+1},x_{l+2},\ldots , x_{p-1}\})=\emptyset. \eqno (32)
$$
Now using (31) and (32), we obtain
$$
d(y)=id(y)+od(y)= \sum^{l-1}_{j=1}
(|E(y\rightarrow x_j)|+|E(x_{j+1}\rightarrow y)|)$$ $$+|E(y\rightarrow x_l)|+|E(x_{1}\rightarrow y)| \leq l+1.
$$
This together with (32) and condition $(M_1)$ implies that $d(x_{l+j})\geq 2p-l-2$ for all $j\in [1,p-l-1]$. Put $P:= x_{1}x_{2}\ldots x_l$. Notice that $|V(P)|=l$ and for every $j\in [1,p-l-1]$, 
$$
 2(p-1)-|V(P)|=2p-l-2\leq d(x_{l+j}).
$$
Therefore, applying Lemma 3.4, we obtain an $(x_{1},x_l)$-path, say $Q$, with vertex set \\ $V(C_{p-1})$. Thus we have
$D(p,3)=[x_{1}yx_l; Q]\subseteq G$. 

Assume next that $G$ contains no cycle of length $p-1$. Then, by Theorem 1.12, either $p$ even and $G$ is isomorphic to $K^*_{p/2,p/2}$ or
$K^*_{p/2,p/2}-\{u\}$, where $u$ is an arbitrary arc of $K^*_{p/2,p/2}$; or $G\in \Phi_p^{p-1}$. Assume that $D$ is isomorphic to $K^*_{p/2,p/2}-\{u\}$, with partite sets $\{x_1,x_2, \ldots ,x_n\}$ and $\{y_1,y_2, \ldots ,y_n\}$ ($p=2n$). Without loss of generality we assume that $u=y_1x_1$. Then $D(p,3)=[x_1y_1x_2; x_1y_2x_3y_3\ldots x_ny_nx_2]$. 

Assume now that  $G\in \Phi_p^{p-1}$. By the definition of $\Phi_p^{p-1}$, $G$ has a Hamiltonian cycle  $x_1x_px_{p-1}\ldots $ $x_2x_1$ such that the vertices $x_1, x_{p-1}$ are not adjacent  and the arcs 
$x_{p-2} x_{p-1}$, $x_{p-1}x_{p}$ are in $E(G)$. Therefore, $D(p,3)=[x_{p-2} x_{p-1}x_{p}; x_{p-2}x_{p-3}\ldots x_2 x_1x_p]$. Thus, in all possible cases  $G$ contains a $D(p,3)$. Theorem 1.19 is proved.   \fbox

\section {Note added in the translation}

 Later on Bang-Jensen, Gutin and Li \cite{[18]}, Manoussakis  \cite{[23]} (see also \cite{[17]}), Bang-Jensen, Guo and Yeo \cite{[19]} proved   the following sufficient conditions (Theorems 5.1-5.4 bellow) for a digraph to be  Hamiltonian. 

\textbf{Theorem 5.1} (Bang-Jensen, Gutin and Li \cite{[18]}). {\it  Let $D$ be a strong digraph of order $n\geq 2$. Suppose that $min\{d(x),d(y)\}\geq n-1$ and  $d(x)+d(y)\geq 2n-1$ for every pair of non-adjacent vertices $x,y$ with a common in-neighbor.  Then $D$ is Hamiltonian.}

\textbf{Theorem 5.2} (Bang-Jensen, Gutin and Li \cite{[18]}). {\it  Let $D$ be a strong digraph of order $n\geq 2$. Suppose that $min\{d^+(x)+d^-(y),d^-(x)+d^+(y)\}\geq n$ for every pair of non-adjacent vertices $x,y$ with a common out-neighbor or a common in-neighbor. Then $D$ is Hamiltonian.}

\textbf{Theorem 5.3} (Bang-Jensen, Guo and Yeo \cite{[19]}). {\it  Let $D$ be a strong digraph of order $n\geq 2$. Suppose that $d(x)+d(y)\geq 2n-1$ and  $min\{d^+(x)+d^-(y),d^-(x)+d^+(y)\}\geq n-1$ for every pair of nonadjacent vertices $x,y$ with a common out-neighbor or a common in-neighbor.  Then $D$ is Hamiltonian.}\\

It is easy to see  that Theorem 5.1 (respectively, Theorem 5.2) implies Ghouila-Houri's (respectively, Woodall's) theorem  and Theorem 5.3 generalizes Theorem 5.2. \\

\textbf{Theorem 5.4} (Manoussakis \cite{[23]}). {\it Let $D$ be a  strong digraph  of order $n\geq 4$. Suppose that $D$ satisfies the following conditions:  For every triple of vertices $x,y,z$ such that $x$ and $y$ are non-adjacent:

 (a) If there is no arc from $x$ to $z$, then $d(x)+d(y)+d^+(x)+d^-(z)\geq 3n-2$. 
 
 (b) If there is no arc from $z$ to $x$, then $d(x)+d(y)+d^-(x)+d^+(z)\geq 3n-2$. 
 
 Then $D$ is Hamiltonian.} \\

Note that Woodall's theorem is an immediate consequence of Theorem 5.4.\\

We pose the following problem: 

\textbf{Problem}. {\it Characterize those  digraphs which satisfy the condition of Theorem 5.1 (or 5.2 - 5.4) but has no Hamiltonian bypass.}\\

In  \cite{[20]}, Darbinyan and Karapetyan the following theorem proved:

\textbf{Theorem 5.5} (Darbinyan and Karapetyan \cite{[20]}).  {\it Let $D$ be a  strong digraph  of order $n\geq 4$. Suppose that $min\{d(x),d(y)\}\geq n-1$ and  $d(x)+d(y)\geq 2n-1$ for every pair of non-adjacent vertices $x,y$ with a common in-neighbor. If  the minimum out-degree of $D$ is at least two and the minimum in-degree of $D$ is  at least three, then $D$ contains a Hamiltonian bypass.}\\

We believe that Theorem 5.5 also is true if we require that the minimum in-degree is at least two instead of three.

\textbf{Theorem 5.6} (Darbinyan \cite{[21]}) . {\it Let $D$ be a  strong digraph  of order $n\geq 4$. Suppose that $D$ satisfies the following conditions:  For every triple of vertices $x,y,z$ such that $x$ and $y$ are non-adjacent: 

(a) If there is no arc from $x$ to $z$, then $d(x)+d(y)+d^+(x)+d^-(z)\geq 3n-2$. 

(b) If there is no arc from $z$ to $x$, then $d(x)+d(y)+d^-(x)+d^+(z)\geq 3n-2$. 

Then $D$ contains a Hamiltonian bypass unless $D$ is isomorphic to the tournament $T_5$.} \\

\textbf{Theorem 5.7} (Darbinyan \cite{[22]}). {\it Let $D$ be a 2-strong digraph  of order $n\geq 3$. Suppose that $d(x)\geq n$   for every  vertex $x\in V(D)\setminus \{x_0\}$, where  $x_0$ is a vertex of $D$. If $D$ is Hamiltonian or $d(x_0)\geq 2(n-1)/5$ then $D$ contains a Hamiltonian bypass.} \\

\end{document}